\documentclass[12pt]{article}
\usepackage{amssymb,amsmath,amsthm,amscd,latexsym}
\usepackage{mathrsfs}
\usepackage{mathrsfs}
\usepackage{amsfonts}
\usepackage{amsmath}
\usepackage{amssymb}
\usepackage{amscd}
\usepackage{mathrsfs}
\usepackage{amssymb}
\usepackage{amsmath}
\usepackage{amsthm}
\usepackage{color}
\usepackage{latexsym}
\usepackage{indentfirst}
\usepackage{enumitem}
\usepackage{anysize}
\usepackage{bbm}
\usepackage[normalem]{ulem}
\usepackage{soul}

\usepackage{cancel}

\renewcommand{\paragraph}{\roman{paragraph}}
\input cyracc.def

\newcommand{\EE}{\mathbb{E}}
\newcommand{\F}{\mathbb{F}}

\newtheorem{theorem}{Theorem}
\newtheorem{corollary}[theorem]{Corollary}

\newtheorem{proposition}[theorem]{Proposition}

\theoremstyle{definition}
\newtheorem{definition}[theorem]{Definition}

\newcommand{\RR}{\mathbb{R}}

\newcommand{\N}{\mathbb{N}}

\begin{document}
\title{\bf Designs in finite metric spaces: a probabilistic approach
\thanks{This research is supported by National Natural Science Foundation of China (12071001, 61672036), Excellent Youth Foundation of Natural Science Foundation of Anhui Province (1808085J20), the Academic Fund for Outstanding Talents in Universities (gxbjZD03).
}}

\author{
\small{Minjia Shi$^{1}$, {Olivier Rioul$^{2}$}, Patrick Sol\'e$^{3,4}$}\\ 
\and 
\small{${}^1$Key Laboratory of Intelligent Computing \& Signal Processing of Ministry of Education,}\\
\small{School of Mathematical Sciences, Anhui University,}\\
\small{Hefei 230601,  China}\\
  \and \small{${}^2$Telecom ParisTech, Palaiseau, France}
\and
\small{${}^3$I2M,(Aix-Marseille Univ., Centrale Marseille, CNRS), Marseille, France}\\
\small{${}^4$ Corresponding Author}
}


\date{}
\maketitle
\begin{abstract} A finite metric space is called here distance degree regular if its distance degree sequence is the same for every vertex. A notion of designs in such spaces is introduced that generalizes
that of designs in $Q$-polynomial distance-regular graphs. An approximation of their cumulative distribution function, based 
on the notion of Christoffel function in approximation theory is given.
As an application we derive limit laws on the weight distributions of binary orthogonal arrays of strength going to infinity. An analogous result
for combinatorial designs of strength going to infinity is given.
\end{abstract}


{\bf AMS Math Sc. Cl. (2010):} Primary 05E35, Secondary O5E20,  05E24 
\section{Introduction}
In a celebrated paper \cite{S} Sidelnikov proved that the weight distribution of binary codes of dual distance
$d^\perp$ going to infinity with the length is close deviates from the normal law up to a term in inverse square root of the dual distance 
\cite[Chap. 9, \S 10]{MS}. Since the times of Delsarte \cite{D}, it is known that the quantity $d^\perp-1$ is the strength of the code viewed as an orthogonal array. Further, this is a special case of designs in so-called $Q$-polynomial association schemes \cite{BH,BCN,D,MS}. The name designs comes
from the Johnson scheme where this notion coincides with that of classical combinatorial designs \cite{D}. Later a similar connection was found between designs in association schemes and designs in lattices \cite{D2,St}. These kinds of generalized designs are popular now in view of the applications in random network coding \cite{E}. In view of this deep connection it is natural to seek to extend
Sidelnikov's theorem to designs in other $Q$-polynomial association schemes than the Hamming scheme. This is a vast research program which might take several years to accomplish.

In the present paper our contribution is twofold. 

Firstly, we develop a theory of designs in finite metric spaces that replaces  the concept of designs in $Q$-polynomial association schemes, when the considered metric space does not afford that structure. We observe that, in contrast with Delsarte 
definition of a design in a $Q$-polynomial association scheme, our definition (Definition 1 below), has an immediate combinatorial meaning in terms of distribution of distances. To wit, the combinatorial meaning of designs in certain $Q$-polynomial association schemes was only derived in \cite{St} in 1986, thirteen years after Delsarte introduced designs in $Q$-polynomial association
schemes in \cite{D}. In particular, the example of permutations with distance the Hamming metric cannot be handled in the context of $Q$-polynomial association schemes, but can  be treated within  our framework. This space had been studied extensively
in the context of permutations codes \cite{BCD,T}. The notion of $t$-design in that space is related to $t$-transitive permutation groups (Theorem 8).

Secondly, we use the technique of Chebyshev-Markov-Stieltjes inequalities in conjunction with orthogonal polynomials to control the difference between the cumulative distribution function of weights 
in designs with that of weights in the whole space. While this technique has been applied by Bannai to the weight distribution of spherical designs \cite{B}, it has not appeared in the literature of algebraic combinatorics so far. 
The bounding quantity in that setting is the Christoffel function, the inverse of the confluent Christoffel-Darboux kernel.
While it is easy to make this quantity explicit in low strength cases, it is difficult to find asymptotic bounds. 
In the present paper, we will use the bounds of Krasikov on the Christoffel-Darboux kernel
of binary Krawtchouk polynomials \cite{K,KL} to derive an alternative proof of the Sidelnikov Theorem. We will give a proof of an analogous result for combinatorial designs by using the limiting behavior of Hahn polynomials.

The material is organized as follows. The next section contains background material on metric spaces, distance-regular graphs and $Q$-polynomial association schemes. Section 3 introduces the definitions that are essential to our approach.
Section 4 contains the main equivalences of our notion of designs in finite metric spaces. Section 5 develops the bounds on the cumulative distribution functions of designs. Section 6 contains some asymptotic results. Section 7 recapitulates the results obtained and gives some significant open problems.
\section{Background material}
\subsection{Metric spaces}
Throughout the paper we write $X$ for a finite set equipped with a metric $d,$ that is to say a map $X\times X \to \N$ verifying the following three axioms
\begin{enumerate}
\item $\forall x,y \in X, \,d(x,y)=d(y,x)$
\item $\forall x,y \in X, \,d(x,y)=0$ iff $x=y$
\item $\forall x,y,z \in X, \,d(x,y)\le d(x,z)+d(z,y).$
\end{enumerate}

In particular if $X$ is the vertex set of a graph the shortest path distance on $X$ is a metric.

The {\em diameter} of a finite metric space is the largest value the distance may take.
A finite metric space is {\em Distance Degree Regular} (DDR) if for every integer $i$ less than the diameter the number $|\{y \in X \mid d(x,y)=i\}|$ is a constant $v_i$ that does not depend of the choice of $x\in X.$

{\bf Example:}
Consider the symmetric group on $n$ letters $S_n$ with metric
$$d_S(\sigma,\theta)=n-F(\sigma \theta^{-1}),$$
where $F(\nu)$ denotes the number of fixed points of $\nu.$
The space $(S_n,d_S)$ is a DDR metric space. If $D_m=m!\sum_{j=1}^m\frac{(-1)^j}{j!}$ denotes the number of fixed-point-free permutations of $S_m,$ ( the so-called d\'erangement number), then $$v_i={ n\choose i}D_{i}.$$
Note that $d_S$ is not a shortest path distance since $ d_S(\sigma,\theta)=1$ is impossible.
Codes in $(S_n,d_S)$ were studied in \cite{T} by using the conjugacy scheme of the group $S_n.$ However, in contrast with the next two subsections, this scheme is neither induced by a graph nor $Q$-polynomial.

\subsection{Distance-regular graphs}
All graphs in this article are finite, undirected, connected, without multiple edges. The neighborhood $\Gamma(x)$ is the set of vertices connected to $x$.
The {\em degree} of a vertex $x$ is the size of $\Gamma(x)$.
A graph is {\em  regular} if every vertex has the same degree. The $i$-neighborhood $\Gamma_i(x)$ is the set of vertices at geodetic distance $i$ to $x$.
The \emph{diameter}  of the graph, denoted by $d$  is the maximum $i$ such that for some vertex $x$ the set  $\Gamma_i(x)$
is nonempty.
A graph is {\em  distance degree regular} (DDR for short) if all graphs $\Gamma_i,$ for $i=1,\dots,d$ are regular.
A graph is {\em  distance regular} (DR for short) if for every two vertices $u$ and $v$ at distance $i$ from each other the values
$
 b_i=| \Gamma_{i+1}(u)\cap \Gamma(v)|$, $
 c_i=| \Gamma_{i-1}(u)\cap \Gamma(v)|$
 depend only on $i$ and do not depend on the choice of $u$ and $v$.
In this case, the graphs $\Gamma_i$ are regular of degree $v_i$ and we will refer to the $v_i$s as the {\em valencies} of $\Gamma$; 
the sequence $\{b_0,\ldots, b_{\mathrm{diam}-1}; c_1,\ldots, c_{\mathrm{diam}}\}$ is usually called the {\em intersection array} of $\Gamma.$
Thus every DR graph is DDR but not conversely.

{\bf Examples}
For background material on the following two examples we refer to \cite{BCN,D,MS}.
\begin{enumerate}
\item The {\em Hamming graph} $H(n,q)$ is a graph on $\F_q^n$ two vertices being connected if they differ in exactly one coordinate. This graph is DR with valencies
$$v_i={n \choose i} (q-1)^i.$$
\item The {\em Johnson graph} $J(\nu,d)$  is a graph on the subsets of cardinality $d$ of a set of cardinality $\nu.$ (Assume $2d<\nu$).
Two subsets are connected iff they intersect in exactly $d-1$ elements.
This graph is DR with valencies
$$v_i={d \choose i}{\nu-d \choose i} .$$
Note that $J(\nu,d)$ can be embedded in $H(\nu,2)$ by identifying subsets and characteristic vectors.
\end{enumerate}
\subsection{$Q$-polynomial association schemes}
An {\em association scheme} on a set $X$ with $s$ classes is a partition of the cartesian product $X\times X=\cup_{i=0}^sR_i$ with the following properties
\begin{enumerate}
\item $R_0=\{(x,x) \mid x\in X\}$
\item $(x,y)\in R_k,$ iff $(y,x)\in R_k,$
\item if $(x,y)\in R_k,$ the number of $z \in X$ such that $(x,z)\in R_i,$ and $(z,y)\in R_j,$ is an integer $p^k_{ij}$ that depends on $i,j,k$ but not on the special choice of $x$ and $y$
\end{enumerate}
A consequence of axiom $3$ is that each graph $R_i$ is regular of degree $v_i,$ say.
It can be shown that the adjacency matrices $D_k$ of the relations $R_k$ span a commutative algebra over the complex with idempotents $J_j$ \cite[Chap. 21]{MS}. Let $\mu_j=$rank$(J_j).$
The {\em first eigenvalues} $p_k(i)$ of the scheme are defined by $D_kJ_i=p_k(i)J_i.$ Considering the matrix $P=(p_k(i))$ and writing $PQ=|X|I,$ with $I$ an identity matrix
defines the {\em second eigenvalues} $q_k(i)$  of the scheme by the relation $Q_{ik}=q_k(i).$
A scheme is said to be $Q$-polynomial if there are numbers $z_0,z_1,\dots,z_s$ such that $q_k(i)=\Phi_k(z_i)$ for some polynomials $\Phi_k$ of degree $k.$
In view of the orthogonality relation \cite[Chap 21, (17)]{MS}
$$\sum_{k=0}^s v_k q_i(k)q_j(k)=|X| \mu_i \delta_{ij},$$ we see that the $\frac{\Phi_i(z)}{\sqrt{\mu_i}}$ form a system of orthonormal polynomials for the scalar product
$$\langle f,g\rangle =\sum_{k=0}^s \frac{v_k}{|X|} f(z_k)g(z_k).$$

{\bf Examples:} In both Hamming and Johnson schemes we have $z_k=k.$
\begin{enumerate}
\item If $\Gamma=H(n,q)$ then $\Phi_k(z)=\frac{K_k(z)}{\sqrt{v_k}},$ where $K_k$ is the Krawtchouk polynomial of degree $k$ given by
the generating function
$$\sum_{k=0}^nK_k(x)z^k=(1+(q-1)x)^{n-x}(1-z)^x.$$
\item If $\Gamma=J(\nu,n)$ then $\Phi_k(z)=\frac{H_k(z)}{\sqrt{v_k}},$ where $H_k$ is the Hahn polynomial of degree $k$ given, as per \cite[(19) p.2481]{DL}, by the formula
$$H_k(z)=m_k\sum_{j=0}^k (-1)^j \frac{{k \choose j}{{\nu+1-k} \choose j} }{{n \choose j}{\nu-n \choose j}} {z \choose j}, $$ where $m_k={ \nu \choose k}-{\nu\choose {k-1}}.$
\end{enumerate}
\section{Preliminaries}
For any integer $\mathcal{N}>0,$ denote by  $[0..\mathcal{N}]$  the set of integers in the range $[0,\mathcal{N}].$
A finite metric space $(X,d)$ is {\em distance degree regular} (DDR) if its distance degree sequence is the same for every point.
Assume $(X,d)$ to be of diameter $n.$ 
In that case $(X,d)$ is DDR iff for each $0\le i\le n$ the graph $\Gamma_i=(X,E_i)$ which connects vertices at distance $i$ in $(X,d)$ is regular of degree $v_i.$ 
Thus $E_0=\{(x,x) \mid x \in X\}$ is the  diagonal of $X^2.$
Note that the $E_i$'s form a partition of $X^2.$

If $D$ is any non void subset of $X$ we define its {\em frequencies} as
$$\forall i \in [0..n], \, f_i=\frac{|D^2\cap E_i|}{|D|^2}.$$ Thus $f_0=\frac{1}{|D|},$ and $\sum\limits_{i=0}^n f_i=1.$ Note also that if $D=X,$ then $f_i=\frac{v_i}{|X|}.$
Consider the random variable $a_D$ defined on $D^2$ with values in $[0..n]$ which affects to an equiprobably chosen 
$(x,y)\in D^2$ the only $i$ such that $(x,y)\in E_i.$ Thus the frequencies $f_i=Prob(a_D=i).$ Denote by $\EE()$ mathematical expectation. Thus
$$\EE(a_D^i)=\sum_{j=0}^n f_jj^i,\,\EE(a_X^i)=\sum_{j=0}^n \frac{v_j}{v}j^i. $$

{\definition The set $D \subseteq X$ is a {\em $t$-design} for some integer $t$ if $$\EE(a_D^i)=\EE(a_X^i)$$ for $i=1,\dots,t.$
}

(Note that trivially $\EE(1)=1$ so that we do not consider $i=0.$)
Thus, distances in $t$-designs are very regularly distributed. For a $2$-design, for instance, the average and variance of the distance coincide with that of the whole space.
 We will see in the next section that in the case of Hamming and Johnson graphs,
we obtain  classical combinatorial objects: block designs, orthogonal arrays.

{\definition We define a scalar product on $\RR[x]$ attached to $D$  by the relation
$$\langle f,g \rangle_D=\sum_{i=0}^n f_i f(i)g(i). $$
Thus, in the special case of $D=X$ we have
$$\langle f,g \rangle_X=\frac{1}{|X|}\sum_{i=0}^n v_i f(i)g(i). $$}
We shall say that a sequence  $\Phi_i(x)$  of  polynomials of degree $i$ is {\em orthonormal of size $N+1$} if it satisfies
$$\forall i,j \in [0..N],\,  \langle \Phi_i,\Phi_j \rangle_X =\delta_{ij},$$
where $N\le n, $ the letter $\delta$ denotes the Kronecker symbol. That sequence is uniquely defined if we assume the leading coefficient of all $\Phi_i(x)$ for $i=0, 1,\dots, N$ to be positive.


For a given DDR metric space $(X,d),$ we shall denote by $N(X)$ the largest possible such $N.$ For instance if $X$ is an $n$-class $P$- and $Q$-polynomial association scheme, it is well-known that $N(X)=n.$
This fact is extended to DDR graphs in the next Proposition.

{\proposition \label{fonda} If none of the $v_i$'s are zero, then $\langle , \rangle_X$ admits an orthonormal system of polynomials of size $n+1.$
In particular, the metric space of a DDR graph admits an orthonormal system of polynomials of size $n+1.$

}

\begin{proof}
 By Lagrange interpolation we see that the functions $1,x,\dots,x^n$ are linearly independent on $[0..n].$ The sequence of the $\Phi_i$'s for $i=0,1,\dots,n$ is then constructed by the usual Gram-Schmidt
 orthogonalization process. Note that this is possible because the bilinear form $\langle , \rangle_X$  is then nondegenerate: $\langle f,f \rangle_X=0 \Rightarrow f=0.$ 
 By properties of the shortest path distance, the property of non vanishing of the  $v_i$'s holds in particular for the metric space of a  DDR graph.
\end{proof}

{\definition For a given $D\subseteq X$ the {\em dual frequencies} are defined for $i=0,1,\dots,N(X)$ as $$\widehat{f_i}=\sum_{k=0}^n\Phi_i(k)f_k.$$}

{\definition For a given $D\subseteq X$ the {\em cumulative distribution function} (c.d.f.) is defined as $$F_D(x)=Prob(a_D\le x)=\sum\limits_{i\le x}f_i.$$}

{\bf Examples:}
\begin{enumerate}
 \item If $D$ is a linear code of $H(n,q),$ with weight distribution $$A_i=\vert \{ x \in D \mid w_H(x)=i\}\vert,$$
then
$F_D(x)=\frac{\sum_{i\le x} A_i}{\vert D\vert}.$
\item If $D$ is a set of points in $J(\nu,k),$ with Hamming distance distribution $B_{2i}$ in $H(\nu,2),$ then 
$F_D(x)=\frac{\sum_{i\le x} B_{2i}}{\vert D\vert}.$

\end{enumerate}

\section{Structure theorems}
First, we give a characterization of $t$-designs in terms of dual frequencies.
{\proposition Let $t$ be an integer $\in [1..N(X)].$ The set $D\subseteq X$ is a $t$-design iff $\widehat{f_i}=0$ for $i=1,\dots,t.$}

\begin{proof}
Note first that $$\EE(a_D^i)=\langle x^i, 1 \rangle_D, \, \EE(a_X^i)=\langle x^i, 1 \rangle_X.$$
Moving the basis of $\RR[x]$ from the $\Phi_i$'s to the basis of monomials we see
that $D$ is a $t$-design if and only if for $i=1,2,\dots,t,$ we have $$ \langle \Phi_i, 1 \rangle_D=\langle \Phi_i, 1 \rangle_X.  $$

Now, by definition, the dual frequency $\widehat{f_i}=\langle \Phi_i, 1 \rangle_D.$
By orthogonality of the $\Phi_i$'s for the scalar product $\langle .,. \rangle_X,$ we see that $\langle \Phi_i,1 \rangle_X=0$ for $i=1,2\dots,t.$
Thus, the condition $\widehat{f_i}=0$ for $i=1,2,\dots,t,$ is equivalent to the fact that $D$ is a $t$-design.

\end{proof}

Next, we connect the notion of designs in $Q$-polynomial association schemes with our notion of designs in metric spaces.
{\theorem If $(X,d)$ is the metric space induced by a $Q$-polynomial DR graph $\Gamma,$ with $z_k=k$ for $k=0,1,\dots,n,$ then 
a $t$-design in $(X,d)$ is exactly a $t$-design in the underlying association scheme of $\Gamma.$}

\begin{proof}
In that situation the frequencies are proportional to the inner distribution (see \cite[p.54]{BCN}) of $D$ in the scheme of the graph, and the dual frequencies are proportional to the dual inner distribution since the second eigenvalues of the scheme, by 
the $Q$-polynomiality condition, are orthogonal polynomials w.r.t. the distribution $\frac{v_i}{|X|}.$ The result follows.
\end{proof}

{\bf Examples:} The following two examples of interpretation of $t$-designs as classical combinatorial objects were observed first in \cite{D} and can be read about in \cite[chap. 21]{MS}.
\begin{enumerate}
\item If $\Gamma$ is the Hamming graph $H(n,q)$ then a $t$-design is an orthogonal array of strength $t.$ That means that every row 
induced by a  $t$-uple columns of $D$ sees the $q^t$ possible values a
constant number of times.
\item If $\Gamma$ is the Johnson graph $J(\nu,n)$ then a $t$-design $D$ is a combinatorial design of strength $t.$ This means the following. Consider $D$ as a collection of subsets of size $n$, traditionally called blocks. That means that every $t$-uple of elements of the groundset is contained in the same number $\eta$ of blocks.
    One says that $D$ is a $t-(\nu,n,\eta)$ design.
\end{enumerate}

Now, we give an example of $t$-design in a metric space that is not a DR graph, or even a DDR graph.

{\theorem If $D\subseteq S_n$ is a $t$-transitive permutation group then it is a $t$-design in $(S_n,d_S).$ }

\begin{proof}
The moments of order $i\le t$ of the number of fixed points of the permutations in $D$ coincide with those of a Poisson law of parameter one. This is a result of Frobenius (1904). A modern exposition is in \cite[Chap. 5.5]{JK}.
The result follows by the definition of $d_S.$
\end{proof}
\section{Distribution functions}
In this section we show that the distribution function of designs are close to that of the whole space.
The proof of the following result follows the philosophy of \cite{B}.
{\theorem Let $D$ be a $t$-design in $(X,d),$ with $t\le N(X).$  Put $\kappa=\lfloor \frac{t}{2}\rfloor.$ Denote by $\lambda(x)$ the Christoffel 
function given by $\lambda(x)=(\sum\limits_{i=0}^\kappa\Phi_i(x)^2)^{-1}.$
Then we have the bound
$$| F_D(x)-F_X(x) |\le \lambda(x).$$
 }
 \begin{proof}
 By Definition 1, we have $$\langle x^i,1 \rangle_X= \langle x^i,1 \rangle_D \, \text{for } \, i=0,1,2,\dots, t.$$
 The orthonormal polynomials for $\langle, \rangle_X$ exist for degrees $\le t$  by the hypothesis $t\le N(X).$ 
 A coincidence of moments up to order $t$ entails a coincidence of orthonormal polynomials up to degree $\kappa$ by Chebyshev determinant for 
 orthonormal polynomials \cite[Lemma 2.1]{H} (see also \cite[(2.2.6), p. 27]{Sz}). By the same formula, the orthonormal polynomials for $\langle, \rangle_D$
 are well-defined for degrees $\le \kappa,$  since 
 the orthonormal polynomials up to degree $\kappa$ attached
 to $X$  exist.
By \cite[Th. 4.1]{H} or \cite[Th. 7.2]{S} we have the Markov-Stieltjes inequalities
\begin{equation}\label{1}
\sum_{x_i<x}\lambda(x_i)\le F_D(x)\le \sum_{x_i<x}\lambda(x_i)+\lambda(x),
\end{equation}
where the $x_i$'s are the $\kappa$ zeros of $\Phi_\kappa(t).$
Similarly we have
\begin{equation}\label{2}
\sum_{x_i<x}\lambda(x_i)\le F_X(x)\le \sum_{x_i<x}\lambda(x_i)+\lambda(x).
\end{equation}
The result follows upon combining equations (\ref{1}) and (\ref{2}).
\end{proof}
 As a bound uniform in $x,$ we have the following result.
 
 {\corollary If $D$ is a binary orthogonal array of strength at least five, then its c.d.f. is close to that of the binomial distribution as
 $$ |F_D(x)-F_X(x) |< \frac{2(n-1)}{3n-2}$$
 }
 \begin{proof}
  We compute explicit lower bounds on $1/\lambda$ by using the first three Krawtchouk polynomials \cite[Chap. 5 \S 7]{MS} given by
 $$K_0=1,K_1(x)=n-2x,K_2(x)=2x^2-2nx+{n \choose2}.$$ We are seeking a lower bound for
 $$1+\frac{K_1(x)^2}{n}+\frac{K_2(x)^2}{{n \choose 2}},$$ when $x\in [0,n].$
 Making the change of variable $y =n-2x \in [-n,n],$ we obtain $$K_2(x)= \frac{y^2-n}{2}, $$ and, therefore
 $$1/\lambda=1+\frac{y^2}{n}+\frac{(y^2-n)^2}{2n(n-1)}=\frac{3n^2-2n+(y^2-1)^2-1}{2n(n-1)},$$
 an increasing function of $y^2$ that takes its minimum over $[0,n^2]$ at $y=0.$
 \end{proof}

{\bf Example:} If $D$ is the extended Hamming code of length $n=16,$  dual distance $8,$ the weight distribution is, in Magma notation \cite{M}, equal to
$$[ <0, 1>, <4, 140>, <6, 448>, <8, 870>, <10, 448>, <12, 140>, <16, 1> ].$$

For $x=8,$ we get $F_D(x)=\frac{1+140+448+870}{2^{11}}\approx 0.712,$ and
$F_X(x)=\frac{\sum\limits_{j=0}^8{n \choose j}}{2^{11}}\approx 0.598.$ The difference is $\approx 0.112 < \frac{3\times 15}{46}\approx 0.652.$

 We give three bounds that are not uniform in $x.$ First, for orthogonal arrays.

 {\corollary If $D$ is a $q$-ary orthogonal array of strength at least two, then its c.d.f. is as close to that of the binomial distribution as
 $$ |F_D(x)-F_X(x) |<\frac{n}{n+(n(q-1)-qx)^2}.$$ }

 \begin{proof}
 Immediate from the data of the first two Krawtchouk polynomials: $$K_0=1,\,K_1(x)=n(q-1)-qx.$$
 \end{proof}
 
 {\bf Example:} If $D$ is a binary Simplex code of length $n=2^m-1,$ there is a unique nonzero weight, namely $\frac{n+1}{2}$ that appears $|D|-1=n$ times.
 If we compute the bound for $x=\frac{n+1}{2}$, its right hand side is $\frac{n}{n+1},$ which is also the value of $F_D(x)$ while 
$F_X(x)> 0.5.$

 Next, we consider combinatorial designs.

 {\corollary If $D$ is a $2-(\nu,n,\lambda)$ design, then its c.d.f. is as close to that of the hypergeometric distribution as
 $$ |F_D(x)-F_X(x) |<\frac{(\nu-n)^3}{(\nu-n)^3+n(\nu-1)^2}.$$  }

 \begin{proof}
 From the data of the first two Hahn polynomials: $$H_0=1,\,H_1(x)=(\nu-1)(1-\frac{\nu x}{n(\nu-n)}),$$
 we obtain $$1/\lambda=1+\frac{H_1(x)^2}{n(\nu-n)}=1+(\nu-1)^2\frac{(n^2-n \nu+\nu x)^2}{n^3(\nu-n)^2},$$
 a monotonic function of $x.$
 \end{proof}

 Eventually, we consider permutation groups. Exceptionally, we do not consider the distance but the codistance $n-d_s.$

 {\corollary If $D$ is a $2$-transitive permutation group on $n$ letters, then the c.d.f. of its fixed points $G_D(x)$ is as close to that of the Poisson law of parameter one
 $P(x)=\sum\limits_{1\le i\le x}\frac{1}{i!}$
  as
 $$ |G_D(x)-P(x) |<\frac{n}{n+(1-x)^2}.$$  }

 \begin{proof}
 Immediate from the data of the first two Charlier polynomials $C_0=1,\, C_1(x)=1-x,$ obtained from the generating series
 $$e^t(1-t)^x=\sum_{n=0}^\infty C_n(x)\frac{t^n}{n!} $$ of \cite[(1.12.11)]{KS}.
 \end{proof}
\section{Asymptotic results}
\subsection{Orthogonal arrays}
In this section we give an alternative proof of a result of Sidelnikov on the weight enumerator of long codes \cite{S}.
 We prepare for the proof by a form of the Central Limit Theorem for the binomial distribution.
 Denote by $\Psi(x)=\frac{1}{\sqrt{2\pi}}\int_{-\infty}^x\exp(-\frac{t^2}{2}) dt$ the cumulative distribution function of the centered normal law of variance unity.
 Let $B_n(x)=\sum_{i\le x}\frac{{n \choose i}}{2^n}$ denote the cumulative distribution function of the binomial law (sum of $n$ Bernoulli trials).
 {\theorem \label{BLT} For some absolute constant $C>0,$ we have
 $$|B_n(x)-\Psi(x)| \le \frac{C}{\sqrt{n}} .$$
  }
 \begin{proof}
 Immediate by Berry-Essen theorem \cite{F2}.
 \end{proof}

 Recall the binary entropy function \cite{MS} defined as
 $$H(x)=-x\log_2 x-(1-x) \log_2(1-x).$$ A tedious but straightforward consequence of Stirling formula is
\begin{equation}\label{entropic}
 {N \choose \alpha N}\sim \frac{2^{NH(\alpha)}}{\sqrt{2\pi \alpha(1-\alpha)}}
\end{equation}
 for $N \to \infty$ and $0<\alpha <1.$ See (1) in \cite{G}.

 {\theorem \label{bigO} Let $n \to \infty,$ and let $k$ be an integer such that $k\sim \theta n,$ with $0<\theta<1$ a real constant. Assume $x=\frac{n}{2}+O(\sqrt{n}).$ Then any binary orthogonal array $D$ with $n$ columns, of strength $\ge 2k+1,$  satisfies
  $$|F_D(x)-B_n(x) |=O(\frac{1}{\sqrt{n}}).$$
  }

 \begin{proof} (sketch)
 We use \cite[Th. 1.1]{K} or \cite[Lemma 4]{KL} to claim the lower bound
 $$\frac{1}{\lambda(x)} =\Omega\big(\frac{{\scriptstyle (k+1)}}{2{n \choose k}}G_k(x)\big),  $$

 where $$G_k(x)=\frac{2k(2x+p_k-n)(n-k)^2\Gamma(x)\Gamma(n-x)}{n^3 (p_k+2)\Gamma(\frac{n}{2}+1)\Gamma(\frac{n}{2}-1)}{n/2 \choose k/2}^2,$$
 where $p_k=2\sqrt{k(n-k)}$ (note that $\mu_k$ tends to a constant in $n.$)

 To derive the said bound, divide numerator and denominator by $n^4$ and simplify. Observe that $p_k\sim 2n\sqrt{\theta(1-\theta)}.$
 For the term ${n/2 \choose k/2}$ we use the entropic estimate mentioned above. We write $\Gamma(x)\Gamma(n-x)=\frac{{(n-2)!}}{{n-2\choose {x-1}}}.$
 We use the Moivre-Laplace formula to get
 $$ \frac{{n \choose x}}{2^n}\sim \frac{\exp(-\frac{(x-n/2)^2}{n/2})}{\sqrt{\pi n/2}}=O(\frac{1}{\sqrt{n}}), $$
 where the constant implied by $O()$ is independent of $x.$
 The result follows after tedious but straightforward manipulations.
 \end{proof}

 We are now ready for the main result of this section.
 {\theorem
 Let $n \to \infty,$ and let $k$ be an integer such that $k\sim \theta n,$ with $0<\theta<1$ a real constant. Assume $x=\frac{n}{2}+O(\sqrt{n}).$ Then any binary orthogonal array $D$ with $n$ columns, of strength $\ge 2k+1,$  satisfies
 $$|F_D(x)-\Psi(x)| =O(\frac{1}{\sqrt{n}}) .$$
  }
 \begin{proof}
 Immediate by combining Theorem \ref{BLT} with Theorem \ref{bigO}.
 \end{proof}
 \subsection{Designs}
Note, before doing asymptotics on the strength of designs, that $t$ designs exist for all $t$ \cite{Te}.
  Let $H_k(x)$ denote the Hahn polynomial of degree $k$ of the variable $x$, as defined in \cite{DL}.
Let $v_k={n \choose k}{{\nu-n} \choose k}$ be the valency of order $k$ of the Johnson graph $J(\nu,n).$
We normalize $\widehat{H_k(x)}=\frac{H_k(x)}{\sqrt{v_k}}.$
{\theorem \label{li} Assume both $\nu$ and $n$ go to infinity with $n/\nu \to p \in (0,1).$ Put $q=1-p.$
Let $z=n x$ with $x \in (0,1).$ Then, we have for fixed $k,$ and $ n\to \infty$ the limit
 $$ \widehat{H_k(z)} \to \frac{(1-x/q)^k}{\sqrt{p^kq^k}}.$$
}

\begin{proof}
First, note that $m_k\sim \frac{\nu^k}{k!}\sim \frac{n^k}{p^k k!}.$ Next, observe that $v_k\sim \frac{n^k (\nu-n)^k}{k!^2}\sim \frac{n^{2k}}{(k!)^2}(q/p)^k.$
This yields $\sqrt{v_k}\sim \frac{n^{k}}{k!}(\sqrt{q/p})^k.$ Combining we obtain $\frac{m_k}{\sqrt{v_k}}\sim 1/\sqrt{p^kq^k}.$
Similar calculations give the term of order $j$ of $H_k(xn)$ to have the limit

$$ (-1)^j {k \choose j}\frac{(\nu^j)n^jx^j}{n^j(\nu-n)^j} \to (-1)^j {k \choose j} (x/q)^j, $$

and, summing on $j$ yield

 $$H_k(xn)\to \sum_{j=0}^k {k\choose j}(-x/q)^j=(1-x/q)^k.$$
 
 The result follows upon writing $ \widehat{H_k(xn)}=\frac{m_k}{\sqrt{v_k}}H_k(xn).$
 \end{proof}

We can now derive the main result of this section.
{\theorem  Let $D$ be a $t-(\nu,n,\eta)$ design with $\nu, n,t\to \infty,$ and $t$ fixed and $ n\sim p\nu$ with $0<p<1$ real constants. Put $q=1-p,$ and $k=\lfloor t/2\rfloor.$
Let $J(\nu,n;x)=\sum_{i \le x} \frac{v_i}{{\nu \choose n}}.$ 
Then
$$|F_D(x)-J(\nu,n;x)|\le \lambda_k(n), $$

where $$\lim_{n \to \infty}{\lambda_k(n)}= \frac{1-a(x)}{1-a(x)^{k+1}}.$$

and $a(x)= \frac{(1-x/q)^2}{\sqrt{pq}} .$}

\begin{proof}
Immediate by taking the limit of the Christoffel-Darboux  kernel of order $k$ given by $$\sum_{j=0}^k{\widehat{H_j(xn)}^2}= \sum_{j=0}^k\frac{H_j(xn)^2}{v_j}$$ and  summing the geometric series
of ratio $a(x)$ coming from Theorem \ref{li}.
 \end{proof}
\section{Conclusion}
In this paper we have used a probabilistic approach to approximate the c.d.f. of designs in various finite metric spaces. 
The key tool is the Christoffel-Darboux kernel attached to the orthonormal
polynomials w.r.t. the valencies of the space. We have used some strong analytic bounds on this quantity for binary Krawtchouk polynomials
derived in \cite{K,KL}. 
It would be desirable to extend these analytical results to other families of polynomials, beginning with $q$-ary Krawtchouk polynomials. 
This special case would yield an alternative proof of the $q$-ary version of Sidelnikov theorem proved by us in
\cite{SRS}. Further, it is a worthwhile project to derive analogous results for the polynomials relevant to the eight types of designs in \cite{St}.
A first step in that direction would be to extend or adapt Theorem 7 to $Q$-polynomial schemes where $z_k$ is not always equal to $k.$
Regarding more general DDR metric spaces, it would be nice to have examples of $t$-designs in the space of permutations that are not $t$-transitive permutation groups.


\begin{thebibliography}{99}
\bibitem{B} E. Bannai, On the weight distribution of spherical $t$-designs, Europ. J. Combinatorics, {\bf 1}, (1980), 19-26.
\bibitem{BH}
A.E.~Brouwer, W.H.~Haemers,
{\em Spectra of graphs,}
Springer (2011).
\bibitem{BCD}I. Blake, G.D. Cohen, M. Deza, Coding with permutations, Information and Control {\bf 43}, (1979), 1--19.
\bibitem{BCN}A.E. Brouwer, A.M. Cohen, A. Neumaier, {\em Distance-regular graphs}, Springer Verlag, Berlin (1989).
\bibitem{D}P. Delsarte, An algebraic approach to the association schemes of Coding Theory, {\it Philips Research Report Suppl. } {\bf 10}, (1973).
\bibitem{D2} P. Delsarte, Association schemes and $t$-designs in regular semi-lattices, J. of Combinatorial Theory A {\bf 20}, (1976), 230--243.
\bibitem{D3} P. Delsarte, Hahn polynomials, discrete harmonics, and $t$-designs, SIAM J. of appl. Math. {\bf 34}, (1), (1978), 157--166.
\bibitem{DL}P. Delsarte, V. Levenshtein, Association schemes and Coding Theory, IEEE Trans. on Information Th., {\bf 44}, ( 6), (1998) 2477--2503.
\bibitem{E} T. Etzion, Problems in $q$-Analogs in Coding Theory, {\tt https://arxiv.org/pdf/1305.6126.pdf}.
\bibitem{F} W. Feller, {\em Introduction to Probability theory and ist Applications}, vol. I, Wiley, New-York (1968).
\bibitem{F2} W. Feller, {\em Introduction to Probability theory and ist Applications}, vol. II, Wiley, New-York (1968).
\bibitem{G}D. Galvin, Three lectures on Entropy and Counting, {\tt https://arxiv.org/pdf/1406.7872.pdf}.
\bibitem{H} W. H\"urliman, An explicit version of the Chebyshev-Markov-Stieltjes inequalities and its applications, J. of Inequalities and Applications (2015).
\bibitem{JK} G. James, A. Kerber, {\it The representation theory of the symmetric group}, Addison-Wesley, Reading MA (1981).
\bibitem{K} I. Krasikov, Bounds for the Christoffel-Darboux kernel of the binary Krawtchouk polynomials, in {\em Codes and Association Schemes}, A. Barg, S. Litsyn eds, AMS (2001), 193--198.
\bibitem{K2} I. Krasikov, Nonnegative Quadratic Forms and Bounds on Orthogonal Polynomials, Journal of Approximation Theory, {\bf 111},(1), (2001), 31--49.
\bibitem{KL} I. Krasikov, S. Litsyn, On the Distance Distributions of BCH
Codes and Their Duals, Designs, Codes and Cryptography, {\bf 23}, (2001), 223--231.
\bibitem{M}{\tt http://magma.maths.usyd.edu.au/magma/}
\bibitem{MS} F.J. MacWilliams,  N.J.A. Sloane, {\it The theory of error-correcting codes}. North-Holland , Amsterdam, (1977).
\bibitem{KS} R. Koekkoek, R.F. Swarttouw, The Askey scheme of hypergeometric orthogonal polynomials and its $q$-analogue, {\tt https://arxiv.org/abs/math/9602214}.
\bibitem{S} V.M. Sidelnikov, Weight spectrum of binary Bose-Chaudhuri-Hocquenghem codes, Problemy Pederachi Informatsii,{\bf 7}, (1), (1971) 14--22.
\bibitem{SRS}M. Shi, O. Rioul, P. Sol\'e, On the asymptotic normality of $Q$-ary linear codes, IEEE Communication Letters {\bf 23}, (11),
(2019), 1895--1898.
\bibitem{St} D. Stanton, $t$-designs in classical association schemes, Graphs and Combinatorics {\bf 2}, (1986) 283--286.
 \bibitem{Sz} G. Szeg\"o, {\it Orthogonal polynomials}, AMS Colloqium Publications {\bf 23}, Providence RI, 4th edition (1975).
 \bibitem{T} H. Tarnanen, Upper bounds on permutation codes via linear programming, Europ. J. of Combinatorics {\bf 20}, (1999), 101--114.
 \bibitem{Te}L. Teirlinck, Non-trivial $t$-designs without repeated blocks exist for all t, Discrete Math.
{\bf 65}, (1987), 301--311
\end{thebibliography}
\end{document}